\documentstyle[12pt,subeqn]{article}
\newfont{\twelvemsb}{msbm10 scaled\magstep1}
\newfont{\eightmsb}{msbm8}
\newfont{\sixmsb}{msbm6}
\newfam\msbfam
\textfont\msbfam=\twelvemsb
\scriptfont\msbfam=\eightmsb
\scriptscriptfont\msbfam=\sixmsb
\catcode`\@=11
\def\Bbb{\ifmmode\let\next\Bbb@\else
  \def\next{\errmessage{Use \string\Bbb\space only in math mode}}\fi\next}
\def\Bbb@#1{{\Bbb@@{#1}}}
\def\Bbb@@#1{\fam\msbfam#1}
\newfont{\twelvegoth}{eufm10 scaled\magstep1}
\newfont{\tengoth}{eufm10}
\newfont{\eightgoth}{eufm8}
\newfont{\sixgoth}{eufm6}
\newfam\gothfam
\textfont\gothfam=\twelvegoth
\scriptfont\gothfam=\eightgoth
\scriptscriptfont\gothfam=\sixgoth
\def\frak{\ifmmode\let\next\frak@\else
  \def\next{\errmessage{Use \string\frak\space only in math mode}}\fi\next}
\def\frak@#1{{\fam\gothfam{{#1}}}}
\catcode`@=12
\newcommand{\be}{\begin{equation}}
\newcommand{\ee}{\end{equation}}
\newcommand{\bea}{\begin{eqnarray}}
\newcommand{\ena}{\end{eqnarray}}
\newcommand{\beo}{\begin{eqnarray*}}
\newcommand{\eno}{\end{eqnarray*}}
\newcommand{\sect}[1]{\setcounter{equation}{0}\section{#1}}

\newcommand{\cle}{\Big[}
\newcommand{\cri}{\Big]}
\newcommand{\ale}{\Big\{}
\newcommand{\ari}{\Big\}}

\newcommand{\RR}{{\Bbb R}}
\newcommand{\ZZ}{{\Bbb Z}}
\newcommand{\ZP}{{\ZZ'}}

\newcommand{\eps}{\varepsilon}

\newcommand{\hs}[1]{\hspace{#1 mm}}
\newcommand{\vs}[1]{\vspace{#1 mm}}
\newcommand{\half}{\frac{1}{2}}
\newcommand{\smbox}[1]{\ \mbox{#1}\ }
\newcommand{\medbox}[1]{\quad\mbox{#1}\quad}
\newcommand{\bigbox}[1]{\qquad\mbox{#1}\qquad}

\newcommand{\cM}{{\cal M}}
\newcommand{\cU}{{\cal U}}
\newcommand{\ta}{{\tilde a}}
\newcommand{\tK}{{\tilde K}}
\newcommand{\wA}{{\widehat A}}
\newcommand{\wC}{{\widehat C}}
\newcommand{\uqa}{{\cU_q(\wA_{N-1})}}
\newcommand{\uqc}{{\cU_q(\wC_N)}}

\newtheorem{theorem}{Theorem}

\begin{document}
\newpage
\pagestyle{empty}
\setcounter{page}{0}
\vfill
\begin{center}
  
  {\LARGE {\bf {\sf ANYONIC REALIZATION OF THE 

                \vs{5}

                QUANTUM AFFINE LIE ALGEBRAS }}}\\[1cm] 
  
  \vs{10}
  
  {\large L. Frappat$^{a}$,
    A. Sciarrino$^{b}$,
    S. Sciuto$^{c}$, 
    P. Sorba$^{a}$
  }

  \vs{10}
  
  {\em $^{a}$ Laboratoire de Physique Th\'eorique ENSLAPP 
       \footnote{URA 1436 du CNRS, associ\'ee \`a l'Ecole Normale 
       Sup\'erieure de Lyon et \`a l'Universit\'e de Savoie.}}
    \\
  {\em Chemin de Bellevue BP 110, F-74941 Annecy-le-Vieux Cedex, France}

  \vs{5}

  {\em $^{b}$ Dipartimento di Scienze Fisiche, Universit\`a di Napoli
       ``Federico II''}
     \\
  {\em and I.N.F.N., Sezione di Napoli, Italy}

  \vs{5}

  {\em $^{c}$ Dipartimento di Fisica Teorica, Universit\`a di Torino}
     \\
  {\em and I.N.F.N., Sezione di Torino, Italy}

\end{center}

\vfill

                                % ABSTRACT
\begin{abstract}
We give a realization of quantum affine Lie algebras $\uqa$ and $\uqc$
in terms of anyons defined on a one-dimensional chain (or on a
two-dimensional lattice), the deformation parameter $q$ being related
to the statistical parameter $\nu$ of the anyons by $q = e^{i\pi\nu}$.
In the limit of the deformation parameter going to one we recover the
Feingold-Frenkel \cite{FF} fermionic construction of undeformed affine
Lie algebras.
\end{abstract}

\vfill
\noindent
Proceedings of the 5th International Colloquium on Quantum Groups and
Integrable Systems, Prague (Czech Republic), 20-22 June 1996, to appear
in the Czechoslovak Journal of Physics, Talk presented by L. Frappat,
\begin{center}
and
\end{center}
Proceedings of the 10th International Conference on Problems of Quantum
Field Theory, Alushta (Crimea, Ukraine), 13-18 May 1996, Talk presented
by S. Sciuto. 

\vfill

\rightline{ENSLAPP-AL-602/96}
\rightline{DSF-T-31/96}
\rightline{DFTT-38/96}
\rightline{q-alg/9607023}
\rightline{July 1996}

\newpage
\pagestyle{plain}

\indent

%SECTION INTRODUCTION
\sect{Introduction}
\label{intro}
\indent

Anyons are particles with {\em any} statistics that interpolate between
fermions and bosons. 

In the first quantization scheme, the notion of statistics is related to
the symmetry properties of the wave function of $N$ identical
particles, the bosons (resp. fermions) corresponding to symmetric
(resp. antisymmetric) wave functions under the exchange of particles.
In two dimensions, there exist more possibilities: the anyons. In that
case, the wave function of $N$ identical particles picks up a {\em
phase factor} under the exchange of particles.
More precisely, if we consider a system of $N$ identical hard-core
particles in $d$ dimensions, the configuration space is $\cM_{N,d} =
[(\RR^d)^N-\Delta]/S_N$ where $\Delta$ is the set of points of
$(\RR^d)^N$ with at least two equal coordinates and $S_N$ is the
permutation group of $N$ elements. The fundamental group 
$\pi_1(\cM_{N,d})=S_N$ for $d > 2$ but $\pi_1(\cM_{N,d})=B_N$ for $d=2$
where $B_N$ is the {\em braid group}. Anyons appear as abelian
representations of the braid group 
\footnote{for a review on anyons see for instance \cite{Lerda}; for a
review on anyonic realization of deformed Lie algebras see
\cite{Varenna}.}.

Let us remark that anyons can consistently be defined also on a
one-dimensional lattice. For simplicity, this will be used here to
construct anyonic realizations of the quantum affine Lie algebras $\uqa$
and $\uqc$. In sect. \ref{conclusion} the construction by means of
two-dimensional anyons \cite{F3S96} will be briefly recalled.

%SECTION A_N
\sect{The unitary case: $\uqa$}
\label{uqa}
\indent

The strategy to construct an anyonic realization of the quantum affine
Lie algebra $\uqa$ with non-vanishing central charge will be the
following: (1) start from the description of $\uqa$ in the
Serre--Chevalley basis, (2) find a fermionic realization of $\wA_{N-1}$
in terms of creation and annihilation operators, (3) construct anyonic
oscillators on a one-dimensional lattice and (4) replace the fermionic
oscillators by suitable anyons in the expressions of the simple
generators of $\uqa$ in the Serre--Chevalley basis.

\medskip

Let us denote by $H_\alpha$ and $E_\alpha^\pm$ where
$\alpha=0,1,\dots,N-1$ the simple generators of $\uqa$ in the
Serre--Chevalley basis. The commutation relations are:
\subequations
\bea
&& \cle H_\alpha,H_\beta \cri = 0 \label{eqA12a} \\
&& \cle H_\alpha,E^\pm_\beta \cri 
= \pm a_{\alpha\beta} E^\pm_\beta \label{eqA12b} \\
&& \cle E^+_\alpha,E^-_\beta \cri 
= \delta_{\alpha\beta} ~ [H_\alpha]_{q_\alpha} \label{eqA12c} 
\ena
\endsubequations
and the quantum Serre relations read as:
\be
\sum_{\ell=0}^{1-a_{\alpha\beta}} (-1)^\ell 
\left[ \begin{array}{c} 1-a_{\alpha\beta} \\ \ell \end{array} 
\right]_{q_\alpha} 
\left( E_\alpha^\pm \right)^{1-a_{\alpha\beta}-\ell} E_\beta^\pm
\left( E_\alpha^\pm \right)^{\ell} = 0
\label{eqA13}
\ee
where the notations are the standard ones, i.e.
\be
[x]_q = \frac{q^x - q^{-x}}{q - q^{-1}} \;, \qquad 
\left[ \begin{array}{c} m \\ n \end{array} \right]_q
= \frac{[m]_q!}{[n]_q![m-n]_q!} \;, \qquad [m]_q! = [1]_q \dots [m]_q
\label{eqA13bis}
\ee
$a_{\alpha\beta}$ being the Cartan matrix of $\wA_{N-1}$ and $q_\alpha=q$
for all $\alpha=0,1,\dots,N-1$.

\medskip

A fermionic realization of $\uqa$ in terms of creation and annihilation
operators is obtained by taking an infinite number of fermionic
oscillators $c_\rho(r), c_\rho^\dagger(r)$ with $\rho=1,\dots,N$ and $r
\in \ZP \equiv \ZZ+\half$, which satisfy the anticommutation relations 
\be
\ale c_\rho(r),c_\sigma(s) \ari 
= \ale c_\rho^\dagger(r),c_\sigma^\dagger(s) \ari = 0
\bigbox{and} \ale c_\rho^\dagger(r),c_\sigma(s) \ari 
= \delta_{\rho\sigma} \delta_{rs}
\label{eqA1}
\ee
the number operator being defined as usual by 
$n_\rho(r) = c_\rho^\dagger(r)c_\rho(r)$. 
\\
These oscillators are equipped with a normal ordered product such that
\be
:c_\rho^\dagger(r) c_\sigma(s): = \left\{ \begin{array}{ll}
c_\rho^\dagger(r) c_\sigma(s) & \smbox{if} s > 0 \\
- c_\sigma(s) c_\rho^\dagger(r) & \smbox{if} s < 0
\end{array} \right.
\label{eqA2}
\ee
and therefore
\be
:n_\rho(r): = \left\{ \begin{array}{ll}
n_\rho(r)  & \smbox{if} r > 0 \\
n_\rho(r) - 1  & \smbox{if} r < 0
\end{array} \right.
\label{eqA2bis}
\ee
Then a fermionic oscillator realization of the simple generators of
$\wA_{N-1}$ in the Serre--Chevalley basis is given by
($\alpha=0,1,\dots,N-1$) (we use small letters for the simple
generators of the undeformed Lie algebra $\wA_{N-1}$)
\be
h_\alpha = \sum_{r \in \ZP}h_\alpha(r) \bigbox{and}
e_\alpha^\pm = \sum_{r \in \ZP} e_\alpha^\pm (r) \label{eqA3a}  
\ee
where $(i=1,\dots,N-1)$ 
\subequations
\bea
&& h_i(r) = n_i(r) - n_{i+1}(r) = ~ :n_i(r): - :n_{i+1}(r):
\label{eqS3a} \\
&& h_0(r) = n_N(r) - n_1(r+1) = ~ :n_N(r): - :n_1(r+1): 
+ \delta_{r,-1/2} \label{eqS3b} \\
&& e_i^+(r) =  c_i^\dagger(r) c_{i+1}(r) \,, \hs{18}
 e_0^+ =  c_N^\dagger(r) c_1(r+1) \label{eqS3c} \\
&& e_i^-(r) =  c_{i+1}^\dagger(r) c_i(r) \,, \hs{18}
 e_0^-(r) =  c_1^\dagger(r+1) c_N(r) \label{eqS3d} 
\ena
\endsubequations
Inserting Eq. (\ref{eqS3b}) into Eq. (\ref{eqA3a}) and taking into
account that the sum over $r$ can be splitted into a sum of two
convergent series only after normal ordering, one can check that 
\be
h_0 = 1 + \sum_{r \in \ZP} :n_N(r): - \sum_{r \in \ZP} :n_1(r): ~ 
= 1 - \sum_{j=1}^{N-1} h_j \;,
\label{eqS4}
\ee
that is the central charge is $\gamma = 1$.

\medskip

We now introduce anyons defined on a one-dimensional lattice:
\be
a_\rho(r) = K_\rho(r) c_\rho(r) \bigbox{and}
{\ta}_\rho(r) ={\tilde  K}_\rho(r) c_\rho(r)
\qquad (1 \le \rho \le N) \label{eqA4}
\ee 
and similar expressions for the conjugated operator $a_\rho^\dagger(r)$
and $\ta_\rho^\dagger(r)$, where the disorder factors $K_\rho(r)$ and
$\tK_\rho(r)$ are expressed as
\be
K_\rho(r) = q^{-\half\sum_{t \in \ZP} \eps(t-r) :n_\rho(t):} 
\bigbox{and}
{\tilde K}_\rho(r) = q^{\half\sum_{t \in \ZP} \eps(t-r) :n_\rho(t):} 
\label{eqA7}
\ee
using the sign function $\eps(t) = |t|/t$ if $t \ne 0$ and $\eps(0) = 0$.
\\
By a direct calculation, one can prove that the $a$-type operators
satisfy the following braiding relations for $r>s$:
\bea
&& \hs{-10} 
a_\rho(r) a_\rho(s) + q^{-1} a_\rho(s) a_\rho(r) = 0 \qquad
a^\dagger_\rho(r) a_\rho(s) + q\ a_\rho(s) a^\dagger_\rho(r) = 0
\nonumber \\
&& \hs{-10} 
a^\dagger_\rho(r) a^\dagger_\rho(s) + q^{-1} a^\dagger_\rho(s)
a^\dagger_\rho(r) = 0 \qquad
a_\rho(r) a^\dagger_\rho(s) + q\ a^\dagger_\rho(s) a_\rho(r) = 0
\label{eqA71}
\ena
and
\bea
&& a_\rho(r) a^\dagger_\rho(r) + a^\dagger_\rho(r) a_\rho(r) = 1 
\nonumber\\
&& a_\rho(r)^2 =  a^\dagger_\rho(r)^2 = 0
\label{eqA72}
\ena
which shows that the operators $a_{\rho}(r), a_{\rho}^\dagger(s)$ are
indeed anyonic oscillators with statistical parameter $\nu$ such that
$q=e^{i\pi\nu}$. The $\ta$-type anyons have the same statistical
parameter $\nu$ but opposite braiding (and  ordering) prescription.

\begin{theorem}
An anyonic realization of the simple generators of the quantum affine
Lie algebra $\uqa$ with central charge $\gamma = 1$
is given by ($\alpha=0,1,\dots,N-1$)
\be
H_\alpha = \sum_{r \in \ZP} H_\alpha(r)
\bigbox{and}
E^\pm_\alpha = \sum_{r \in \ZP} E^\pm_\alpha(r)
\label{eqA8}
\ee
where $(1 \le i \le N-1)$
\subequations
\bea
&& \hs{-5} H_i(r) = h_i(r) = :n_i(r): - :n_{i+1}(r): \label{eqA10a} \\
&& \hs{-5} E^+_i(r) = a^\dagger_i(r) a_{i+1}(r) \hs{28}
E^-_i(r) = {\ta}^\dagger_{i+1}(r) {\ta}_i(r) \label{eqA10b} \\
&& \hs{-5} H_0(r) = h_0(r) = :n_N(r): - :n_1(r+1): + \delta_{r,-1/2}
\label{eqA10c} \\ 
&& \hs{-5} E^+_0(r) = q^{-\half \eps(r+1/2)} ~ a^\dagger_N(r) a_1(r+1)
\hs{5} 
E^-_0(r) = q^{-\half \eps(r+1/2)} ~ {\ta}^\dagger_1(r+1) {\ta}_N(r) 
\nonumber \\
&& \hs{-5} \label{eqA10d}
\ena
\endsubequations
\end{theorem}

\noindent
{\bf Proof}
We must check that the realization Eqs. (\ref{eqA10a}-\ref{eqA10d}) 
indeed satisfy the quantum affine Lie algebra $\uqa$ in the
Serre--Chevalley basis.
\\
Fist of all, inserting eqs. (\ref{eqA4}), the expressions
(\ref{eqA10a}-d) become
\be
E_\alpha^\pm(r) = e_\alpha^\pm(r) q^{\half\sum_{t\in\ZP} \eps(t-r)
:h_\alpha(t):} \label{eqA9}
\ee
where the generators $h_\alpha(r)$ and $e_\alpha^\pm(r)$ coincide with
those defined in eqs. (\ref{eqS3a}-d) corresponding to the undeformed
affine Lie algebra $\wA_{N-1}$.

Now, consider the {\em non-extended} Dynkin diagram of $\wA_{N-1}$, to
which corresponds the set of generators
$\{h_\alpha(r),e_\alpha^\pm(r)\}$ where $\alpha \ne 0$. For a fixed $r
\in \ZP$, the set $\{h_\alpha(r),e_\alpha^\pm(r)\}$ for $\alpha \ne 0$
is a representation of $A_{N-1}$ of spin 0 and 1/2, and thus of
$\cU_q(A_{N-1})$ \cite{FMS}. Thanks to eqs. (\ref{eqA8}) and
(\ref{eqA9}), $H_\alpha,E_\alpha^\pm$ are the correct coproduct in
$\cU_q(A_{N-1})$. Therefore, the relations (\ref{eqA12a}-c) hold for
$\alpha,\beta \ne 0$ (step 1).

We consider then the {\em extended} Dynkin diagram of $\wA_{N-1}$ and we
delete a dot $\mu$ which is not the affine dot. For a fixed $r \in
\ZP$, the set $\{h_\alpha(r),e_\alpha^\pm(r)\}$ for $\alpha \ne \mu$ is
a representation of $A_{N-1}$ of spin 0 and 1/2, and thus of
$\cU_q(A_{N-1})$. The eqs. (\ref{eqA8}) and (\ref{eqA9}) give once
again the correct coproduct in $\cU_q(A_{N-1})$ and hence
$H_\alpha,E_\alpha^\pm$ form a representation of $\cU_q(A_{N-1})$.
Therefore, the relations (\ref{eqA12a}-c) hold for $\alpha,\beta \ne
\mu$ (step 2).

Finally, in the case $\cU_q(\wA_1)$ the previous arguments fail (in
particular to prove the quantum Serre relations). Such equations can
however be checked explicitly by using the braiding properties of the
anyonic oscillators $a$ and $\ta$.

%SECTION C_N
\sect{The symplectic case: $\uqc$}
\label{uqc}

\indent

Let us consider now the $\uqc$ case. As in the $\uqa$ case, a guideline
is given by the fermionic realization of the (undeformed) affine
algebra. Such a fermionic realization in terms of creation and
annihilation operators is easily obtained by noticing that the folding
of the Dynkin diagram of $\wA_{2N}$ leads to the Dynkin diagram of
$\wC_N$ (see figure below).
%%%%%%%%%%%%%%%%%%%%%%%%%%%%%%%%%%%%%%%%%%%%%%%%%%%%%%%%%%%%
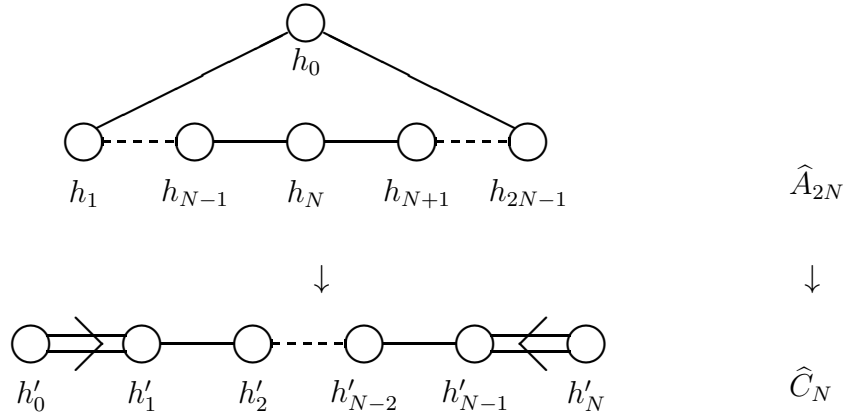
\beo
&
\begin{picture}(180,80)
\thicklines
\multiput(0,20)(42,0){5}{\circle{14}}
\put(84,65){\circle{14}}
\put(7,20){\dashbox{3}(28,0)}
\put(49,20){\line(1,0){28}}
\put(91,20){\line(1,0){28}}
\put(5,25){\line(2,1){72}}
\put(163,25){\line(-2,1){72}}
\put(133,20){\dashbox{3}(28,0)}
\put(0,0){\makebox(0.4,0.6){$h_1$}}
\put(42,0){\makebox(0.4,0.6){$h_{N-1}$}}
\put(84,0){\makebox(0.4,0.6){$h_N$}}
\put(126,0){\makebox(0.4,0.6){$h_{N+1}$}}
\put(168,0){\makebox(0.4,0.6){$h_{2N-1}$}}
\put(84,50){\makebox(0.4,0.6){$h_0$}}
\end{picture}
& \hs{20} \wA_{2N} \\
\\
&\downarrow& \hs{22} \downarrow \\
&& \\
&
\begin{picture}(220,20)
\thicklines
\multiput(0,20)(42,0){6}{\circle{14}}
\put(6,17){\line(1,0){30}}
\put(6,23){\line(1,0){30}}
\put(49,20){\line(1,0){28}}
\put(91,20){\dashbox{3}(28,0)}
\put(133,20){\line(1,0){28}}
\put(174,17){\line(1,0){30}}
\put(174,23){\line(1,0){30}}
\put(27,20){\line(-1,-1){10}}\put(27,20){\line(-1,1){10}}
\put(185,20){\line(1,1){10}}\put(185,20){\line(1,-1){10}}
\put(0,0){\makebox(0.4,0.6){$h'_0$}}
\put(42,0){\makebox(0.4,0.6){$h'_1$}}
\put(84,0){\makebox(0.4,0.6){$h'_2$}}
\put(126,0){\makebox(0.4,0.6){$h'_{N-2}$}}
\put(168,0){\makebox(0.4,0.6){$h'_{N-1}$}}
\put(210,0){\makebox(0.4,0.6){$h'_N$}}
\end{picture}
& \hs{20} \wC_N
\eno
%%%%%%%%%%%%%%%%%%%%%%%%%%%%%%%%%%%%%%%%%%%%%%%%%%%%%%%%%%%%
Denoting the simple generators of $\wA_{2N}$ by $h_i,
e_i^\pm$ with $i=0,\dots,2N-1$, the lifting of the automorphism
associated to the symmetry of the Dynkin diagram of $\wA_{2N}$ leads
to the following expression of the simple generators 
of $\wC_N$:
\bea
&& h'_0 = h_0, \quad h'_i = h_i + h_{2N-i}, \quad h'_N = h_N \nonumber \\
&& {e'}_0^\pm = e_0^\pm, \quad {e'}_i^\pm = e_i^\pm + e_{2N-i}^\pm, \quad
{e'}_N^\pm = e_N^\pm
\label{eqC1}
\ena
In the undeformed case, it is immediate to see that the Serre
relations of $\wC_N$ are satisfied (actually one can write
them entirely in terms of the Serre relations of $\wA_{2N}$).

Now, we go to the deformed case. The main idea of the construction is
to use the folding $\wA_{2N} \longrightarrow \wC_N$ to obtain a 
realization of $\uqc$ in terms of anyons, following with some
modifications the procedure used in Ref. \cite{FMS} to build $U_q(C_N)$.
Using the fermionic oscillators $c_\mu(r), c^\dagger_\mu(r)$ of the
previous Section, with the same normal ordering prescription Eq.
(\ref{eqA2bis}), one defines the following anyons: 
\be
a_\mu(r) = K_\mu(r) c_\mu(r) \bigbox{and}
{\ta}_\mu(r) = {\tilde K}_\mu(r) c_\mu(r) \qquad (1 \le \mu \le 2N)
\label{eqC2}
\ee
where the disorder operators are now
\bea
&& K_\mu(r) = {\tilde K}_\mu^\dagger(r) = q^{ -\half \sum_{t \in \ZP} 
\eps(t-r) (:n_\mu(t): - :n_{2N-\mu+1}(t):) - :n_{2N-\mu+1}(r):} \nonumber
\\ 
&& \hs{50} \medbox{for} \mu = 1,\dots,N \nonumber \\
&& K_\mu(r) = {\tilde K}_\mu^\dagger(r) = q^{ -\half \sum_{t \in \ZP} 
\eps(t-r) (:n_\mu(t): - :n_{2N-\mu+1}(t):) + :n_{2N-\mu+1}(r):} \nonumber
\\ 
&& \hs{50} \medbox{for} \mu = N+1,\dots,2N \label{eqC3}
\ena
\begin{theorem}
An anyonic realization of the simple generators of the quantum affine
Lie algebra $\uqc$ with central charge equal to 1 is given by 
($\alpha = 0,1,\dots,N$) 
\be
H_\alpha = \sum_{r \in \ZP} H_\alpha(r)
\bigbox{and}
E^\pm_\alpha = \sum_{r \in \ZP} E^\pm_\alpha(r)
\label{eqC7}
\ee
where
\bea
&& H_i(r) = h_i(r) = \ :n_i(r): - :n_{i+1}(r): +
:n_{2N-i}(r): - :n_{2N-i+1}(r): \nonumber \\
&& E^+_i(r) = a^\dagger_i(r) a_{i+1}(r) + a^\dagger_{2N-i}(r)
a_{2N-i+1}(r) \nonumber \\
&& E^-_i(r) = {\ta}^\dagger_{i+1}(r) {\ta}_i(r) + 
{\ta}^\dagger_{2N-i+1}(r) {\ta}_{2N-i}(r)
\label{eqC4}
\ena
are associated to the simple short roots of $C_N$ ($i=1,\dots,N-1$),
\bea
&& H_N(r) = h_N(r) = \ :n_N(r): - :n_{N+1}(r): \nonumber \\
&& E^+_N(r) = a^\dagger_N(r) a_{N+1}(r) \nonumber \\
&& E^-_N(r) = {\ta}^\dagger_{N+1}(r) {\ta}_N(r) 
\label{eqC5}
\ena
are associated to the simple long root of $C_N$, and
\bea
&& H_0(r) = h_0(r) = \ :n_{2N}(r): - :n_1(r+1): 
+ \delta_{r,-1/2} \nonumber \\ 
&& E^+_0(r) = a^\dagger_{2N}(r) a_1(r+1) \nonumber \\
&& E^-_0(r) = {\ta}^\dagger_1(r+1) {\ta}_{2N}(r)
\label{eqC6}
\ena
are associated to the affine root of $C_N$.
\end{theorem}

\noindent
{\bf Proof}
We have now to check that the realization Eqs. (\ref{eqC4}),
(\ref{eqC5}), (\ref{eqC6}) indeed satisfy the quantum affine Lie
algebra $\uqc$ in the Serre--Chevalley basis, i.e. the equations
(\ref{eqA12a}-\ref{eqA12c}) for $\alpha,\beta = 0,1,\dots,N$, together
with the quantum Serre relations (\ref{eqA13}) hold, where
$a_{\alpha\beta}$ is now the Cartan matrix of $\wC_N$ and $q_i = q$ for
$i = 1,\dots,N-1$ (short roots) and $q_0 = q_N= q^2$ (long roots).

The proof is based on the strategy exposed in the previous section.
First of all, let us define for $\mu=1,\dots,2N-1$
\subequations
\bea
&& k_\mu(r) = n_\mu(r) - n_{\mu+1}(r) \label{eqC10a} \\
&& f_\mu^+(r) = c_\mu^\dagger(r) c_{\mu+1}(r) \label{eqC10b} \\
&& f_\mu^-(r) = c_{\mu+1}^\dagger(r) c_\mu(r) = \left( f_\mu^+(r)
\right)^\dagger \label{eqC10c}
\ena
and
\bea
&& f_0^+(r) = c_{2N}^\dagger(r) c_1(r+1) \label{eqC10d} \\
&& f_0^-(r) = c_1^\dagger(r+1) c_{2N}(r) = \left( f_0^+(r)
\right)^\dagger \label{eqC10e}
\ena
\endsubequations
As before, by using the definitions (\ref{eqC2}) of $a,\ta$ and
taking into account Eq. (\ref{eqC3}), the expressions (\ref{eqC4}), 
(\ref{eqC5}) and (\ref{eqC6}) simplify as ($\alpha=0,1,\dots,N$)
\be
E_\alpha^\pm(r) = e_\alpha^\pm(r) \ q_\alpha^{\half \sum_{t \in \ZP}
\eps(t-r) h_\alpha(r)} \label{eqC6bis}
\ee
where 
\subequations
\bea
&& \hs{-10} e_i^\pm(r) = f_i^\pm(r) q^{\half k_{2N-i}(r)}
+ f_{2N-i}^\pm(r) q^{-\half k_i(r)} \medbox{for} i=1,\dots,N-1
\label{eqC20a} \\
&& \hs{-10} e_N^\pm(r) = f_N^\pm(r) \label{eqC20b} \\
&& \hs{-10} e_0^\pm(r) = f_0^\pm(r) \label{eqC20c}
\ena
\endsubequations

The step 1 is achieved by noticing that for any fixed $r \in \ZP$, the
set $\{k_i(r),f_i^\pm(r);$ $i=1,\dots,N-1\}$ and
$\{k_{2N-i}(r),f_{2N-i}^\pm(r) ; i=1,\dots,N-1\}$ are representations
of spin 0 and 1/2 of $U_q(A_{N-1})$; therefore the set 
$\{h_i(r),e_i^\pm(r);$ $i=1,\dots,N-1\}$ is a representation of 
$U_q(A_{N-1})$ corresponding to the coproduct of these
two representations. Analogously, $\{h_N(r),e_N^\pm(r)\}$ is a
representation of $A_1$ of spin 0 and 1/2 and thus of $U_q(A_1)$.
Moreover, from the properties of the fermionic oscillators, one can
directly check the quantum Serre relations and Eqs. (\ref{eqA12b}-c)
involving $h_N(r)$ and $e_N^\pm(r)$ and thus prove that the set 
$\{ h_\rho(r), e_\rho^\pm(r);\ \rho = 1,\dots,N \}$ is a representation
of $U_q(C_N)$. Then Eqs. (\ref{eqC7}) and (\ref{eqC6bis}) show that also
the set $\{ H_\rho, E_\rho^\pm;\ \rho = 1,\dots,N \}$ is a
representation of $U_q(C_N)$, obtained by iterated coproduct 
\footnote{Let us remark that this anyonic representation of $U_q(C_N)$
does not require any Gutzwiller projection on the fermionic oscillators
and thus it is an improvement of the one given in ref. \cite{FMS}.}.

For the step 2, one has to check that the equations
(\ref{eqA12a}-\ref{eqA12c}) and (\ref{eqA13}) also hold when
$\alpha,\beta$ can take the value $0$. Actually the only non trivial
cases are for Eqs. (\ref{eqA12c}) and (\ref{eqA13}) when
$\alpha=0,\beta=1$ or $\alpha=1,\beta=0$. Thus we have just to show
that the set $\{ H_\alpha, E_\alpha^\pm;\ \alpha = 0,1 \}$ is a
representation of $U_q(C_2)$; as $h_0(r)$ involves oscillators both of
the site $r$ and of the site $r+1$, it is convenient to rewrite Eq.
(\ref{eqC7}) rearranging the (convergent) series for $H_1$ and
$E_1^\pm$: 
\be
H_1 = \sum_{r \in \ZP} H'_1(r) \qquad
E_1^\pm = \sum_{r \in \ZP} {E'}_1^\pm(r) 
\label{eqC11}
\ee
where
\be
H'_1(r) = h'_1(r) = k_1(r+1) + k_{2N-1}(r) 
\label{eqC20}
\ee
and
\bea
&& {E'}_1^+ = a_1^\dagger(r+1) a_2(r+1) + a_{2N-1}^\dagger(r) a_{2N}(r)
\label{eqC21a} \\
&& {E'}_1^- = \ta_2^\dagger(r+1) \ta_1(r+1) + \ta_{2N}^\dagger(r) 
\ta_{2N-1}(r) \label{eqC21b} 
\ena
Then using the definitions (\ref{eqC2}) of $a,\ta$ and taking into 
account Eq. (\ref{eqC3}), one gets
\be
{E'}_1^\pm(r) = {e'}_1^\pm(r) \ q^{\half \sum_{t \in \ZP} \eps(t-r)
h'_1(t)} \label{eqC14}
\ee
where
\be
{e'}_1^\pm(r) = f_{2N-1}^\pm(r) q^{\half k_1(r+1)} 
+ f_1^\pm(r+1) q^{-\half k_{2N-1}(r)}  
\label{eqC15}
\ee
Let us show now that the set $\{ h_0(r), e_0^\pm(r), h'_1(r),
{e'}_1^\pm(r) \}$ is a representation of $U_q(C_2)$. For a fixed $r \in
\ZP$, the triple $\{ h'_1(r), {e'}_1^\pm(r) \}$ is a representation of
$U_q(A_1)$, corresponding to the coproduct of 
$\{ k_{2N-1}(r),f_{2N-1}^\pm(r) \}$ and $\{ k_1(r+1),f_1^\pm(r+1) \}$.
Therefore one has $[{e'}_1^+(r),{e'}_1^-(r)] = [h'_1(r)]_q$. Using
elementary properties of fermio\-nic oscillators, one checks the
equations $[e_0^+,{e'}_1^-] = [e_0^-,{e'}_1^+] = 0$ and the quantum
Serre relations on the generators $e_0^\pm,{e'}_1^\pm$. Then from Eq.
(\ref{eqC6bis}) for $\alpha=0$ and Eq. (\ref{eqC14}), it follows that
the set $\{ H_0, E_0^\pm, H_1, E_1^\pm \}$ is a representation of
$U_q(C_2)$ and thus Eqs. (\ref{eqA12a}-c) and the Serre relations 
(\ref{eqA13}) on the generators (\ref{eqC7}) also hold for 
$\alpha,\beta \in \{0,1\}$, which concludes the proof. 
As in the $\uqa$ case, the central charge is 1 because of the normal 
ordering prescription Eq.
(\ref{eqA2bis}).

%SECTION CONCLUSION
\section{Conclusion}
\label{conclusion}

\indent

We have presented here a method to get an anyonic realization of the
quantum affine Lie algebras $\uqa$ and $\uqc$ with central charge
$\gamma=1$. Let us emphasize the role of the definition of the ordering
of anyons which is crucial in this construction. Representations with
vanishing central could be built in the same way by using alternative
ordering prescriptions. These representations with $\gamma=0$ and
$\gamma=1$ can be combined together to get representations with arbitrary
positive integer central charges which are in general reducible. 
This composition of representations, with the correct deformed 
coproduct, is naturally achieved by considering anyons defined on a 
two-dimensional lattice \cite{F3S96}. 
It is worthwhile to notice that these anyonic
constructions have nothing to do with $q$-deformed oscillators.
Finally, let us mention that it is also possible to obtain a
supersymmetric version of the construction (see ref. \cite{nous} for
the case $\cU_q(\wA(m,n))$).

\section*{Acknowledgements}

\indent

This research was partially supported by MURST and EU, within the
framework of the program "Gauge Theories, Applied Supersymmetry and
Quantum Gravity" under contract SCI*-CT92-0789. Two of us (L.F. and 
P.S.) are also grateful to EU for partial financial help (European
Network ERBCHRXTCXT 920069).

\end{document}